\documentclass[11pt,oneside,reqno]{amsart}
\usepackage{amsmath,amssymb,amsbsy,amsfonts,amsthm,
latexsym,amsopn,amstext, amsxtra,euscript,amscd,hyperref}
\usepackage{xcolor}
\newtheorem{theorem}{Theorem}[section]
\newtheorem{definition}{Definition}

\newtheorem{corollary}{Corollary}[section]

\begin{document}

\title[]{On the transcendence  of growth constants associated with polynomial recursions}
\author{Veekesh Kumar}

\pagenumbering{arabic}

\begin{abstract} 
Let $P(x):=a_d x^d+\cdots+a_0\in\mathbb{Q}[x]$, $a_d>0$, be a polynomial of degree $d\geq 2$.  Let $(x_n)$ be a sequence of integers satisfying 
\begin{equation*}
x_{n+1}=P(x_n)~~\mbox{for all~~} n=0,1,2\ldots,\quad\mbox{and~~~} x_n\to\infty\quad\mbox{as~~~} n\to\infty.
\end{equation*}

Set $\alpha:=\lim_{n\to\infty} x^{d^{-n}}_n$. Then, under the assumption $a_d^{1/(d-1)}\in\mathbb{Q}$, in a recent result by Dubickas \cite{dubickas}, either $\alpha$ is transcendental, or $\alpha$ can be an integer, or a quadratic Pisot unit with $\alpha^{-1}$ being its conjugate over $\mathbb{Q}$. In this paper, we study the nature of such $\alpha$ without the assumption that $a_d^{1/(d-1)}$ is in $\mathbb{Q}$, and we prove that either the number $\alpha$ is transcendental, or $\alpha^h$ is a Pisot number with $h$ being the order of the torsion subgroup of the Galois closure of the  number field $\mathbb{Q}(\alpha, a_d^{-\frac{1}{d-1}})$.
Other results presented in this paper investigate the solutions of the inequality $||q_1 \alpha_1^n+\cdots+q_k \alpha_k^n +\beta||<\theta^n$ in $(n,q_1,\ldots,q_k)\in \mathbb{N}\times(K^\times)^k$, considering whether $\beta$ is rational or irrational. Here, $K$ represents a number field, and $\theta\in (0,1)$. The notation $||x||$ denotes the distance between $x$ and its nearest integer in $\mathbb{Z}$.
\end{abstract}
\address[Veekesh Kumar]{Department of Mathematics, Indian Institute of Technology, Dharwad  580011, Karnataka, India.}

\email[Veekesh Kumar]{veekeshk@iitdh.ac.in}


\maketitle

\section{Introduction}

\medskip
Let $P(x):=a_d x^d+\cdots+a_0\in\mathbb{Q}[x]$, $a_d>0$, be a polynomial of degree $d\geq 2$.  Let $(x_n)$ be a sequence of integers satisfying 
\begin{equation*}\label{eq1}
\tag{1.1}
x_{n+1}=P(x_n)~~\mbox{for all~~} n=0,1,2\ldots,\quad\mbox{and~~~} x_n\to\infty\quad\mbox{as~~~} n\to\infty.
\end{equation*}

Recently Wagner and Ziegler \cite{wagner} showed that for the above sequence $(x_n)$, $\lim_{n\to\infty} x^{d^{-n}}_n=\alpha$ exists. Moreover, they showed that  $\alpha>1$ and it is either irrational or an integer.  Also, it was shown that  such sequence $(x_n)_n$ takes the form 
\begin{equation*}\label{eq2}
\tag{1.2}
x_n=a_d^{-1/(d-1)}\alpha^{d^n}-\frac{a_{d-1}}{d a_d}+O(\alpha^{-d^n}).
\end{equation*}
The proof of the asymptotic formula \eqref{eq2} has already been discussed in \cite{wagner} and \cite{dubickas}. Here, we provide a sketch of the proof for completeness. By the substitution $y_n:=a_d^{1/d-1}(x_n+\frac{a_{d-1}}{d a_d})$, the recursion becomes 
\begin{align*}
y_{n+1}&=a_d^{1/d-1}\left(P(x_n)+\frac{a_{d-1}}{d a_d}\right)\\
&=a_d^{d/d-1}\left(x_n+\frac{a_{d-1}}{d a_d}\right)+O(x^{d-2}_n)\\
&=y^d_n+O(y^{d-2}_n)\quad \mbox{as}~ n\to\infty.
\end{align*}
Since $x_n\to\infty$ as $n\to \infty$,  so $y_n\to\infty$. By increasing $n$, if necessary we can assume that the sequence $(y_n)_{n\geq 0}$ is increasing and none of the $y_n$ is zero. Taking the logarithm, we get 
\begin{equation*}\label{eq3}
\tag{1.3}
\log \left(\frac{y_{n+1}}{y^d_n}\right)=O(y^{-2}_n) \quad \mbox{as}~ n\to\infty.
\end{equation*}
Express $\log y_n$ as follows: 
$$
\log y_n=d\log y_{n-1}+\log\left(\frac{y_n}{y^d_{n-1}}\right)=d^n\log y_0+\sum_{k=0}^{n-1}d^{n-k-1}\log\left(\frac{y_{k+1}}{y^d_k}\right).
$$
Since the series $\displaystyle \sum_{k=0}^\infty d^{-k-1}\log\left(\frac{y_{k+1}}{y^d_k}\right)$is convergent, we can re-write $\log y_n$ as follows: 
$$
\log y_n=d^n\left(\log y_0+\sum_{k=0}^\infty d^{-k-1}\log\left(\frac{y_{k+1}}{y^d_k}\right)\right)-\sum_{k=n}^\infty d^{n-k-1}\log\left(\frac{y_{k+1}}{y^d_k}\right).
$$
Set 
$$
\log \alpha=\log y_0+\sum_{k=0}^\infty d^{-k-1}\log\left(\frac{y_{k+1}}{y^d_k}\right).
$$
Then, we get 
$$
\log y_n=d^n\log \alpha-\sum_{k=n}^\infty d^{n-k-1}\log\left(\frac{y_{k+1}}{y^d_k}\right). 
$$
Since the sequence $(y_n)n$ is increasing from some point onwards, we have $y_n \leq y{n+1} \leq \cdots$ for sufficiently large $n$. Together with equation \eqref{eq3}, we deduce that
$$
\log y_n=d^n \log \alpha+O\left(\frac{1}{y^2_n}\right),
$$ 
which in turns gives  $y_n=\alpha^{d^n}+O(\alpha^{-d^n})$ as $n\to \infty$, and hence
$$
x_n=a^{-1/d-1}_d \alpha^{d^n}-\frac{a_{d-1}}{d a_d}+O\left(\frac{1}{\alpha^{d^n}}\right) \quad
\mbox{as}~~ n\to\infty.
$$
This completes the proof of \eqref{eq2}.
\smallskip

Very recently, Dubickas \cite{dubickas} studied the transcendence of numbers $\alpha$ under the assumption $a^{1/d-1}_d \in \mathbb{Q}$. More precisely, he showed the possibility of such  $\alpha$  can be an integer, a quadratic Pisot unit with $\alpha^{-1}$ being its conjugate over $\mathbb{Q}$, or a transcendental number. As a consequence, he established the transcendence of several constants given by the polynomial recursion.  For example,  he considered the sequence $1, 2,5, 26, 277, 458330,\ldots$, given by $x_0=1$ and
$$
x_{n+1}=x^2_n+1\quad \mbox{for~~~} n=0,1,2,\ldots.
$$
It can also defined as $x_n=[\kappa^{2^n}]$, $n=0,1,\ldots$, where 
$$
\kappa:=\lim_{n\to\infty} x^{2^{-n}}_n=\prod_{n=1}^\infty\left(1+\frac{1}{x^2_n}\right)^{\frac{1}{2^{n+1}}}.
$$
The sequence $2,3,8,63, 3968,15745023,\ldots$ given by $x_0=2$ and 
$$
x_{n+1}=x^2_n-1\quad \mbox{for}~~ n=0,1,2,\ldots.
$$
Then $x_n=[\zeta^{2^n}]$, $n=0,1,2,\ldots,$ where 
$$
\zeta:=\lim_{n\to\infty} x^{2^{-n}}_n.
$$
Another example is Sylvester’s sequence $2, 3, 7, 43, 1807, 3263443,\ldots$, where $x_0=2$ and 
$$
x_{n+1}=x^2_n-x_n+1\quad \mbox{for~~~} n=0,1,2,\ldots.
$$
In this case $x_n$ can also given by $x_n=[\gamma^{2^n}]$, where $\gamma:=\lim_{n\to\infty} x^{2^{-n}}_n$.
These above sequences can be found in the On-Line Encyclopedia of Integer Sequences \cite{online}(see also \cite{aho}).  By the above mentioned result of Dubickas, the constants $\kappa, \zeta$ and $\gamma$ are transcendental. 
\smallskip

Notice that in both of the above sequences, the respective polynomials have a leading coefficient of $1$, so the condition $a_d^{1/(d-1)}\in\mathbb{Q}$ is satisfied. In this paper, our main result provides a variant of  Dubickas' result without assuming $a_d^{1/(d-1)}\in\mathbb{Q}$. Many fascinating examples of sequences satisfying \eqref{eq1} can be found in Finch's book on mathematical constants \cite{finch}.
\smallskip

Here is our main result:

\begin{theorem}\label{maintheorem4}
 Suppose that an integer sequence $(x_n)_n$ satisfies a recursion of the form $x_{n+1}=P(x_n)$
for some polynomial $P=a_d X^d+a_{d-1}x^{d-1}+\cdots+a_0\in\mathbb{Q}[X]$ degree $d\geq 2$ and $a_d>0$. Assume further that $x_n\to\infty$ as $n\to\infty$.
Then either the number 
$$
\alpha=\lim_{n\to\infty} (x_n)^{\frac{1}{d^n}}
$$
is transcendental, or $\alpha^h$ is a Pisot number, with $h$ being the order of the torsion subgroup of the Galois closure of the  number field $K=\mathbb{Q}(\alpha, a_d^{-\frac{1}{d-1}})$.
Moreover, if $a_d$ is an integer, then either $\alpha$ is transcendental, or $(a_d)^{\frac{d-2}{d-1}}\alpha^{d^m}$ is a Pisot number for some non-negative integer $m$.
\end{theorem}

\noindent{\bf Remark 1.}~  In the proof of the  above theorem, we are using the fact that $d\cdot a_d \cdot a_d^{-1/(d-1)}\alpha^{d^n}$ is pseudo-Pisot number for infinitely many positive integers $n$, which in turn implies that $\alpha$ is an algebraic integer, as proven by either the result of Corvaja and Zannier \cite{corv} or Theorem \ref{maintheorem1} below. So, if we additionally assume that $a_d ^{-1/d-1}\in\mathbb{Q}$ (which is an essential assumption in Dubickas's result), then we can conclude that $\alpha^{d^m}$ is a Pisot number for some integer $m\geq 0$. Without this assumption, the conclusion is no longer true; we can only assert that $\alpha^D$ is a Pisot number for some positive integer $D$,  but this $D$ may not be of the form $d^m$. In our Theorem \ref{maintheorem4}, we obtain a similar conclusion as Dubickas's result, except for the extra factor $a_d^{d-2/d-1}.$
\bigskip

We illustrate Theorem \ref{maintheorem4} with some examples.
\bigskip

\noindent{\bf Example 1.~} Consider the sequence given by $x_1=3$ and $x_{n+1}=x_n^2-2$. It has been shown in \cite[Page 2]{wagner} that 
$$
x_n=L_{2^n}=\left(\frac{1+\sqrt{5}}{2}\right)^{2^n}+\left(\frac{1-\sqrt{5}}{2}\right)^{2^n}
$$
for all $n\geq 1$, where $L_n$ is the $n$th Lucas number.  Thus, the limit of the sequence $x^{2^{-n}}_n$ would be the golden ratio in this case. Since $a_d=1$, by taking $m=0$, we see that $a_d^{(d-2)/(d-1)}\alpha^{d^m}$ is nothing but the golden ratio. Therefore, it is a Pisot number.
\bigskip

\noindent{\bf Example 2.~} Consider the polynomial $P(x) = 2x^d$. Then, the $n$th
term of the sequence defined by $x_0 = 1$ and $x_{n+1} = 2x_n^d$ for $n = 0, 1, 2,...,$ and it is  equal to 
$$
x_n=2^{(d^n-1)/(d-1)}.
$$
Hence, $\alpha=\displaystyle \lim_{n\to\infty} x_n^{d^{-n}}=2^{1/(d-1)}$, which is an algebraic integer. Since the order of the torsion subgroup of the Galois closure of the number field $\mathbb{Q}(\alpha)$ is $d-1$, we have $\alpha^{d-1}=2$, making it a Pisot number. This also explains why the exponent $h$ in Theorem \ref{maintheorem4} is the best possible. Furthermore, since $a_d$ is a positive integer, by taking   $m=1$,  we see that
$$
a_d^{(d-2)/(d-1)}\alpha^{d^m}=2^{(d-2)/(d-1)} 2^{d/(d-1)}=2^2=4,
$$
which is again a Pisot number. 
\bigskip

For a complex number $x$,  $||x||$  denotes the distance of $x$ from its nearest integer in $\mathbb{Z}$.  In other words, 
$$
||x||:=\mbox{min}\{|x-m|:m\in\mathbb{Z}\}.
$$
We recall  the following definition. 
\begin{definition}
A tuple $(\alpha_1, \ldots, \alpha_k)$ of non-zero algebraic numbers is called  {\it non-degenerate} if $\alpha_i/\alpha_j$ is not a  root of unity for all integers $1\leq i < j \leq k$.
\end{definition}
\bigskip

When working with sums of the form $q_1 \alpha^n_1+\cdots+q_k \alpha^n_k$, we can assume that $(\alpha_1,\ldots,\alpha_k)$ is non-degenerate without loss of generality.  To see this, suppose $\frac{\alpha_k}{\alpha_{k-1}}=\zeta$ is an $h$-th root of unity. For $0\leq a\leq h-1$, we restrict to $n\in\mathbb{N}$ congruent to $a$ modulo $h$, and write $n=a+hm$. Then the sum $q_1\alpha^n_1+\cdots+q_k \alpha^n_k$ is equal to the sum $q_1 \alpha^a_1(\alpha^h_1)^{m}+\cdots+(q_{k-1+\zeta^a q_k})\alpha^a_{k-1}(\alpha^{h}_{k-1})^m$, which has fewer terms than the original sum.  
\bigskip

We also recall  the following definition introduced in \cite{kul}.
\begin{definition}\label{def2}
Let $(\beta_1,\ldots,\beta_k)$ be  a tuple of distinct non-zero algebraic numbers. Set 
$$
B:=\{\beta\in\bar{\mathbb{Q}}^\times\backslash\{\beta_1,\ldots,\beta_k\}:~\beta=\sigma(\beta_i)~~\mbox{for some~~}\sigma:\mathbb{Q}(\beta_1,\ldots,\beta_k)\to \mathbb{C}~~\mbox{and~~} 1\leq i\leq k\}.
$$
 Then the tuple $(\beta_1,\ldots,\beta_k)$ is  called pseudo-Pisot if $\sum_{i=1}^k\beta_i+\sum_{\beta\in B}\beta\in\mathbb{Z}$  and $|\beta|<1$  for every $\beta\in B$. Moreover, if $\beta_i$ is an algebraic integer for $1\leq i\leq k$ then the tuple $(\beta_1,\ldots,\beta_k)$ is called Pisot.
\end{definition}
Let $h(x)$ denote the absolute logarithmic Weil height.   By  ${\it sublinear~ function}$, we mean a function $f:\mathbb{N}\rightarrow (0,\infty)$ satisfying $\displaystyle\lim_{n\to\infty}\frac{f(n)}{n}=0$. Let $G_\mathbb{Q}$ be the absolute Galois group of $\mathbb{Q}$.
\bigskip

We require  the following diophantine approximation result of Kulkarni, Mavraki, and Nguyen \cite{kul}, which extends a seminal work by Corvaja and Zannier \cite{corv}. 
\begin{theorem}\label{maintheorem1}{\rm(Kulkarni, Mavraki and Nguyen)}~
Let $r\in\mathbb{N}$, let $(\delta_1,\ldots,\delta_r)$ be a non-degenerate tuple of algebraic numbers  with $|\delta_i|\ge 1$ for $1\leq i\leq r$.  Let $K$ be a number field and  $f$ be a sublinear function.  Suppose for some $\theta\in (0,1)$, the set $\mathcal{M}$ of  tuple  $(n, q_1,\ldots,q_r)\in \mathbb{N}\times(K^\times)^r$ satisfying the inequality 
\begin{equation*}
||q_1 \delta_1^n+\ldots+q_r\delta_r^n||<\theta^n\quad \mbox{and}~~~ \max_{1\leq i\leq k}h(q_i)<f(n)
\end{equation*}
is infinite. Then   the following holds: 
\begin{enumerate}
\item[(i)] $\delta_i$ is an algebraic integer for $i=1,\ldots,r$.
\item[(ii)] For each $\sigma\in G_\mathbb{Q}$  and $1\leq i\leq r$ such that $\frac{\sigma(\delta_i)}{\delta_j}$ is not a root of unity for $1\leq j\leq r$, we have $|\sigma(\delta_i)|<1$.
\end{enumerate}
Moreover for all but finitely many  tuples $(n,q_1,\ldots,q_r)\in\mathcal{A}$
\begin{enumerate}
\item[(iii)]
 $(q_1\delta^n_1,\ldots,q_r\delta^n_r)$ is pseudo-Pisot.
\item[(iv)]
$\sigma(q\delta^n_i)=q_j\delta^n_j$ precisely for  those triples $(\sigma, i, j)\in G_\mathbb{Q}\times\{1,\ldots,r\}^2$ such that $\frac{\sigma(\delta_i)}{\delta_j}$ is  a root of unity. 
\end{enumerate}
\end{theorem}
\vspace{.2cm}

Theorem \ref{maintheorem1} plays a crucial role in proving all the results presented in this paper.

\section{some other results}
Here is our second theorem, which is an immediate consequence of Theorem \ref{maintheorem1}. 
\begin{theorem}\label{maintheorem2}
Let $k\in\mathbb{N}$, and let $(\alpha_1,\ldots,\alpha_k)$ be a non-degenerate tuple of algebraic numbers   with $|\alpha_i|\geq 1$ for $1\leq i\leq k$ and none of $\alpha_i$ is  root of unity. Let $\beta\in (0,1)$ be an algebraic irrational  number. Let $K$ be a number field and  $f$ be a sublinear function.  Then for  any  ~$\theta\in(0,1)$, there are only finitely many tuples $(n,q_1,\ldots,q_k)\in\mathbb{N}\times (K^\times)^k$ satisfying 
\begin{equation*}\label{eq1.1}
\tag{2.1}
||q_1 \alpha_1^n+\ldots+q_k\alpha_k^n+\beta||<\theta^n\quad \mbox{and}~~~ \max_{1\leq i\leq k}h(q_i)<f(n).
\end{equation*}
\end{theorem} 
Very recently, the case $k=1$  with $q_1$ is a fixed algebraic number, independently also proved by the author \cite{kumar}.
We have the following  corollary  of Theorem \ref{maintheorem2}.
\begin{corollary}\label{cor} 
Let $\alpha_1,\ldots,\alpha_k$ be multiplicatively independent  algebraic numbers with $|\alpha_i|\geq 1$ for $1\leq i\leq k$ and none of $\alpha_i$ is a root of unity. Let $P(x_1,\ldots,x_k)$ be a non-zero polynomial with algebraic coefficients and  constant term is irrational.  Suppose  that  for some $\theta\in (0,1)$,  there are infinitely many $n\in\mathbb{N}$ such that $||P(\alpha_1^n,\ldots,\alpha^n_k)||<\theta^n$.  Then,  at least one  of $\alpha_i$ is transcendental.
\end{corollary}

It is natural to ask what can we say in the case when  $\beta$ is a  rational number in Theorem \ref{maintheorem2}. First, let us  see some remarks, then we will come back to this case.
\bigskip

\noindent{\bf Remark 2.1.}
{\it We observe that in the case    when  $\alpha_i$ is pseudo-Pisot number,  $q_i=1$ for all $1\leq i\leq k$, and  $\beta$ is a  real number in the interval  $(0,1)$,  the inequality  \eqref{eq1.1} can have only finitely many solutions in $n$  for any given $\theta\in(0,1)$. This  can be seen as follows:  suppose there are infinitely many  integers $n\geq 1$ and some $\theta'\in(0,1)$ such that the inequality 
\begin{equation*}\label{eq1.2}
\tag{2.2}
||\alpha^n_1+\cdots+\alpha^n_k+\beta||<\theta'^n
\end{equation*}
holds. Let $p_n $ be the nearest integer to $\alpha^n_1+\cdots+\alpha^n_k+\beta$. Then $p_n$ is of the form 
$$
p_n=\mathrm{Tr}_{\mathbb{Q}(\alpha_1)/\mathbb{Q}}(\alpha^n_1)+\cdots+\mathrm{Tr}_{\mathbb{Q}(\alpha_k)/\mathbb{Q}}(\alpha^n_k)+[\beta]+a,
$$
for all sufficiently large values of $n$,  where $a$ is either $0$ or $1$. Here we used the fact that when $\alpha$ is pseudo-Pisot number then the $\mbox{Tr}_{\mathbb{Q}(\alpha)/\mathbb{Q}}(\alpha^n)$  is the nearest integer of  $\alpha^n$ for all sufficiently large positive integers $n$.  Thus we have
$$
||\alpha^n_1+\cdots+\alpha^n_k+\beta||=\left|\sum_{i=2}^{d_1} \alpha^n_{1,i}+\cdots+\sum_{i=2}^{d_k} \alpha^n_{k,i}+\{\beta\}-a\right|.
$$
Since $\{\beta\}\in (0,1)$,  we have  $\{\beta\}-a$ is non-zero. On the other hand,  by the hypothesis $\alpha_i$'s are pseudo-Pisot numbers,    we get that  $\alpha^n_{j,i}\to 0$ for every pair $(i,j)$.  Thus by these observations, we have
\begin{equation*}\label{eq1.3}
\tag{2.3}
||\alpha^n_1+\cdots+\alpha^n_k+\beta||>c(\beta)>0
\end{equation*}
holds for all large positive integers $n$.  From \eqref{eq1.2} and \eqref{eq1.3},  we get a contradiction and hence the assertion. }
\bigskip

\noindent{\bf Remark 2.2.} {\it If we take $q_1 =\frac{1}{2}$, $\alpha_1=\frac{1+\sqrt{5}}{2}$,  $\beta=\frac{1}{2}$ and $\theta=|1-\sqrt{5}/2|$. It can be easily seen that $\mathrm{Tr}(\alpha^{2^n}_1)$ is an odd integer for all $n\in\mathbb{N}$, where  $\mathrm{Tr}:=\mathrm{Tr}_{\mathbb{Q}(\sqrt{5})/\mathbb{Q}}$. Hence, we get that  
$$
\Vert q_1 \alpha^{2^n}_1+\beta\Vert\leq \left|\frac{\alpha^{2^n}_1}{2}-\left(\frac{\mathrm{Tr}(\alpha^{2^n}_1)}{2}+\frac{1}{2}\right)+\frac{1}{2}\right|=\frac{1}{2}\theta^{2^n}
$$
for all sufficiently large values of $n$.   This explains that in general, the assumption $\beta$ is an irrational number cannot be removed in Theorem \ref{maintheorem2}.} 
\bigskip

\noindent{\bf Remark 2.3.} {\it The assumption that none of $\alpha_i$ is the root of unity is a necessary condition in Theorem \ref{maintheorem2}. Take $\alpha_1$ to be any root of unity of order $h$, $\alpha_2=\frac{1+\sqrt{5}}{2}$ and  $q_2=\frac{1}{2}$, then take any $q_1$ and $\beta$ such that the $q_1\alpha_1 + \beta = 1/2$. Then for any $n$ such that $n\equiv 1~\mbox{mod}~ h$ and $\mathrm{Tr}_{\mathbb{Q}(\sqrt{5})/\mathbb{Q}}(\alpha_2^n)$ is odd, we  get back to the situation as in  Remark 2.2.}
\bigskip

In the case when $\beta$ is  rational number, we have the following result.
\begin{theorem}\label{maintheorem3}
Let $k\in\mathbb{N}$, and let $(\alpha_1,\ldots,\alpha_k)$ be a non-degenerate tuple of algebraic numbers   with $|\alpha_i|\geq 1$ for $1\leq i\leq k$ and none of $\alpha_i$ is  root of unity. Let $\beta$ be a non-integral  rational  number. Let $K$ be a number field, $\mathcal{O}_K$ be its ring of integers and  $f$ be a sublinear function.  Then for  any  ~$\theta\in(0,1)$, there are only finitely many  tuples  $(n, q_1,\ldots,q_k)\in \mathbb{N}\times(\mathcal{O}_K^\times)^k$ satisfying the inequality 
\begin{equation*}
||q_1 \alpha_1^n+\ldots+q_k\alpha_k^n+\beta||<\theta^n\quad \mbox{and}~~~ \max_{1\leq i\leq k}h(q_i)<f(n)
\end{equation*}
\end{theorem}
 Remark 2.2 explains our restriction $(n, q_1,\ldots,q_k)\in \mathbb{N}\times(\mathcal{O}_K^\times)^k$ in Theorem \ref{maintheorem3}.

\section{Proof of Theorem \ref{maintheorem4}.}
 As we have discussed in the beginning, the sequence $x_n$ as stated in the theorem can be given by the following asymptotic formula
\begin{equation*}
x_n=a_d^{-1/(d-1)}\alpha^{d^n}-\frac{a_{d-1}}{d a_d}+O(\alpha^{-d^n}),
\end{equation*}
where   $\alpha=\lim_{n\to\infty} x^{\frac{1}{d^n}}_n$ and it  is strictly greater than 1. 
\smallskip

Suppose that $\alpha$ is an algebraic number. Let $L=\mathbb{Q}(a_d^{-1/(d-1)}, \alpha)$, and $K$ be its Galois closure. Let $h$ be the order of the torsion subgroup of $K^\times$. Since $\alpha>1$, we have $ d a_d a_d^{-1/(d-1)} \alpha^{d^n}>1$ for all large enough integers $n$. Then by the hypothesis, the inequality, 
$$
| d a_d\cdot a_d^{-1/(d-1)} \alpha^{d^n}-(d a_d x_n-a_{d-1})|<C(\alpha)\left(\frac{1}{\alpha}\right)^{d^n}
$$
has infinitely many solutions in $n$ for some constant $C(\alpha)>0$.  Since $a_d$ and $a_{d-1}$ are fixed rational numbers, we can multiply by a suitable fixed positive integer to reduce the case when the term $(d a_d x_n - a_{d-1})$ becomes an integer. Therefore, the above inequality can be rewritten as:
\begin{equation*}\label{eq1.5}
\tag{3.1}
\Vert d a_d\cdot a_d^{-1/(d-1)} \alpha^{d^n}\Vert<C(\alpha)\left(\frac{1}{\alpha}\right)^{d^n},
\end{equation*}
which  has infinitely many solutions in $n$ for some constant $C(\alpha)>0$. We denote $\mathcal{A}$ be a set of positive integers $n$ for which  the inequality \eqref{eq1.5} holds. Since $\mathcal{A}$ is infinite, there exists an integer $a \in \{0, 1, \ldots, h-1\}$ such that $d^n = a+hm$ for infinitely many natural numbers $m$. 
\smallskip

By taking $k=1$ with $q_1=d \cdot a_d\cdot (a_d^{-1/(d-1)}) \alpha^a$ and $\delta_1=\alpha^h$, we can observe that the hypothesis of Theorem \ref{maintheorem1} is satisfied. Therefore, according to part (i) of Theorem \ref{maintheorem1}, $\alpha^h$ is an algebraic integer. Additionally, based on part (iii), for all sufficiently large values of $n$ satisfying \eqref{eq1.5} and $d^n=a+hm$, we have $d \cdot a_d\cdot (a_d^{-1/(d-1)}) \cdot \alpha^{d^n}$ is a pseudo-Pisot number. In order to complete the proof of this theorem, it suffices to show that $|\sigma(\alpha^h)|<1$ for each embedding $\sigma\neq \text{Id}:\mathbb{Q}(\alpha^h)\to\mathbb{C}$.

We first observe that any conjugate $\sigma(\alpha^h)\neq \alpha^h$ has an absolute value less than or equal to $1$. Assume that $|\sigma(\alpha^h)|>1$. Since $d \cdot a_d\cdot (a_d^{-1/(d-1)}) \cdot \alpha^{d^n}$ is a pseudo-Pisot number, we must have $\rho(q_1\alpha^{hm})=\rho(q_1)\sigma(\alpha^h)^m=q_1 \alpha^{hm}$ for all but finitely many values of $m$ such that $d^n=a+hm\in\mathcal{A}$ and some $\rho\in\mbox{Gal}(L/\mathbb{Q})$, where $\sigma$ is the restriction of the automorphism $\rho$ to $\mathbb{Q}(\alpha^h)$. Then by property (iv) of Theorem \ref{maintheorem1}, we have $\sigma(\alpha^h)/\alpha^h$ is a root of unity, and hence $\sigma(\alpha^h)=\alpha^h$, which is a contradiction. Thus, we conclude that $|\sigma(\alpha^r)|\leq 1$.

Now, we show that the possibility $|\sigma(\alpha^h)|=1$ cannot occur. If we have $|\sigma(\alpha^h)|=1$, then the quotient $\frac{\sigma(\alpha^h)}{\alpha^h}$ is not a root of unity. By property (ii) of Theorem \ref{maintheorem1}, we have $|\sigma(\alpha^h)|<1$, which contradicts the assumption that $|\sigma(\alpha^h)|>1$. This proves that $|\sigma(\alpha^h)|<1$ for each embedding $\sigma\neq \text{Id}:\mathbb{Q}(\alpha^h)\to\mathbb{C}$, and hence finishes the proof of the first part.
\smallskip

For the moreover part of the theorem, we apply Theorem \ref{maintheorem1} with the inputs $k=1$, $q_1=d \cdot a_d\cdot (a_d^{-1})^{1/(d-1)}$, and $\delta_1=\alpha$. Consequently, we can conclude that $\alpha$ is an algebraic integer, and $d a_d^{(d-2)/(d-1)}\alpha^{d^n}$ is a pseudo-Pisot number for all but finitely many values of $n\in\mathcal{A}$. Since $a_d\in\mathbb{N}$, we can deduce that $a_d^{(d-2)/(d-1)}$ is an algebraic integer. This implies that $a_d^{(d-2)/(d-1)}\alpha^{d^n}$ is also an algebraic integer.

Therefore, given the fact that $d a_d^{(d-2)/(d-1)}\alpha^{d^n}$ is a pseudo-Pisot number for all but finitely many values of $n\in\mathcal{A}$, we can further conclude that $a_d^{(d-2)/(d-1)}\alpha^{d^m}$ is a Pisot number for some non-negative integer $m$.
\section{Proofs of Theorems \ref{maintheorem2}, \ref{maintheorem3} and Corollary \ref{cor}}
\noindent{\bf Proof of Theorem \ref{maintheorem2}.}~The proof of this theorem follows from Theorem \ref{maintheorem1}.  Suppose there exists an infinite set $\mathcal{A}$ of tuples $(n, q_1,\ldots,q_k)\in\mathbb{N}\times (K^\times)^k$, where $K$ is some number field and some $\theta\in (0,1)$ such that the inequality
\begin{equation*}
0<||q_1 \alpha_1^n+\ldots+q_k\alpha_k^n+\beta||<\theta^n\quad \mbox{and}~~~ \max_{1\leq i\leq k}h(q_i)<f(n)
\end{equation*}
holds. Since  none of  $\alpha_i$ is root of unity, we have the tuple $(\alpha_1,\ldots,\alpha_k,1)$ is non-degenerate.   Then we apply  Theorem \ref{maintheorem1} with the inputs  $\mathcal{M}=\mathcal{A}$,  $r=k+1$, $\delta_i=\alpha_i$ for $1\leq i\leq k$  and $q_{k+1}=\beta, \delta_{k+1}=1$, and conclude by part $(iv)$ of Theorem \ref{maintheorem1} that $\beta$ is fixed by all $\sigma\in G_\mathbb{Q}$, and hence $\beta$ is rational. This is a contradiction to the assumption that $\beta$ is an irrational number and hence  the theorem.
\vspace{.3cm}

\noindent{\bf Proof of Theorem \ref{maintheorem3}.}
For the proof of this theorem, we argue  by contradiction. Suppose there exists an infinite set $\mathcal{A}$ of tuples $(n, q_1,\ldots,q_k)\in\mathbb{N}\times (K^\times)^k$, where $K$ is some number field  and some $\theta\in (0,1)$ such that the inequality
\begin{equation*}
0<||q_1 \alpha_1^n+\ldots+q_k\alpha_k^n+\beta||<\theta^n\quad \mbox{and}~~~ \max_{1\leq i\leq k}h(q_i)<f(n)
\end{equation*}
holds. Since  none of $\alpha_i$ is  root of unity, we have the tuple $(\alpha_1,\ldots,\alpha_k,1)$ is non-degenerate.   Then we apply  Theorem \ref{maintheorem1} with the inputs  $\mathcal{M}=\mathcal{A}$, $r=k+1$, $\delta_i=\alpha_i$ for $1\leq i\leq k$  and $q_{k+1}=\beta, \delta_{k+1}=1$, and conclude by part $(i)$ of Theorem \ref{maintheorem1} that each $\alpha_i$ is an algebraic integer. By part $(iii)$ of Theorem \ref{maintheorem1}, together with Defintion \ref{def2}, we also have 
$$
\sum_{i=1}^k\mathrm{Tr}_{\mathbb{Q}(q_i\alpha^n_i)/\mathbb{Q}}(q_i\alpha^n_i)+\beta\in\mathbb{Z}
$$
for infinitely many values of $n\in\mathcal{A}$. Since each $\alpha_i$ is an algebraic integer and by the hypothesis each $q_i$ is also  an algebraic integer, we get that  $q_i\alpha_i$ is an algebraic integer for $i=1,\ldots,k$.  Thus the sum $\sum_{i=1}^k\mathrm{Tr}_{\mathbb{Q}(q_i\alpha^n_i)/\mathbb{Q}}(q_i\alpha^n_i)\in\mathbb{Z}$, which in turns implies that $\beta$ is an integer. This contradicts  the assumption that $\beta$ is not an integer and hence the theorem..
\bigskip

\noindent{\bf{Proof of Corollary \ref{cor}}.}
Let $P(x_1,\ldots,x_k)=\displaystyle\sum_{\bar{i}=(i_1,\ldots,i_k)}a_{\bar{i}}x^{i_1}_1\cdots x^{i_k}_k+a_{\bar{0}}$ be a polynomial with real algebraic coefficients and degree in each variable $x_i$ is $d_i$, and the constant term $a_{\bar{0}}$ is irrational. By the hypothesis  
$$
||P(\alpha_1^n,\ldots,\alpha^n_k)||=||\sum_{\bar{i}}a_{\bar{i}}\left(\alpha^{i_1}_1\cdots\alpha^{i_k}_{k}\right)^n+a_{\bar{0}}||<\theta^n
$$
holds for infinitely many values of $n$.  Since $|\alpha_i|\geq 1$ for all $1\leq i\leq k$ and $\alpha_1,\ldots,\alpha_k$ are multiplicatively  independent, we have the tuple $(\alpha^{i_1}_1\cdots\alpha^{i_k}_k:1\leq i_j\leq d_j, 1\leq j\leq k )$  is non-degenerate. Then by  applying  Theorem \ref{maintheorem2} with inputs $k=d_1 d_2\cdots d_k$,  $\delta_{i_1,\ldots,i_k}=\alpha^{i_1}_1\cdots\alpha^{i_k}_k, \beta=a_{\bar{0}}$ and $q_{i_1,\ldots,i_k}=a_{\bar{i}}$ for $1\leq i_j\leq d_j, 1\leq j\leq k$,  we conclude that at least one of $\alpha_i$ is transcendental.
\bigskip

\noindent{\bf Acknowledgements.}  I am grateful to  the anonymous referee whose constructive suggestions and comments helped in improving the exposition. My sincere thanks to Prof. Dang Khoa Nguyen for several helpful comments and suggestions on an earlier version of this paper, especially for bringing my attention to shortening the proof of Theorem \ref{maintheorem2}. I  am very grateful to Prof. Purusottam Rath for the fruitful discussion and for suggesting some remarks. This work was started during the author's visit to Chennai Mathematical Institute in July 2022, I am thankful to this institute for providing me with an excellent research environment.   Also, I would like to thank Dr. Siddhi S. Pathak for careful reading of an earlier draft of this paper.

\end{document}